\documentclass[12pt]{article}

\usepackage{amsfonts,amsbsy,graphicx}

\newcommand{\ih}{\'{\i}}
\newcommand{\eh}{\hspace{.06in}}
\newcommand{\A}{{\cal{A}}}
\newcommand{\B}{{\cal{B}}}
\newcommand{\cC}{{\cal{C}}}

\newcommand{\K}{{\cal{K}}}

\newcommand{\X}{{\cal{X}}}
\newcommand{\C}{\mathbb{C}}
\newcommand{\E}{\mathbb{E}}
\newcommand{\N}{\mathbb{N}}
\newcommand{\R}{\mathbb{R}}
\newcommand{\Z}{\mathbb{Z}}
\newcommand{\hC}{\hat{\C}}
\newcommand{\m}{\frac{1}{2}}
\newcommand{\af}{\alpha}

\newcommand{\eps}{\varepsilon}
\newcommand{\deh}{\partial}
\newcommand{\ds}{\displaystyle}
\newcommand{\ovl}{\overline}
\newcommand{\BE}{\begin{equation}}
\newcommand{\EE}{\end{equation}}
\newcommand{\fns}{\footnotesize}
\newcommand{\Lim}[1]{\lower5.2pt\hbox{${{\ds\lim}\atop ^{#1}}$}}

\begin{document}
\centerline{\Large{\bf Scherk Saddle Towers of Genus Two in $\R^3$}}
\bigskip

\centerline{M{\fns \'ARCIO} F{\fns ABIANO DA} S{\fns ILVA} \& V{\fns AL\'ERIO} R{\fns AMOS} B{\fns ATISTA}}
\bigskip
\begin{abstract}
In 1996 M. Traizet obtained singly periodic minimal surfaces with Scherk ends of arbitrary genus by desingularizing a set of vertical planes at their intersections. However, in Traizet's work it is not allowed that three or more planes intersect at the same line. In our paper, by a {\it saddle-tower} we call the desingularization of such ``forbidden'' planes into an embedded singly periodic minimal surface. We give explicit examples of genus two and discuss some advances regarding this problem. Moreover, our examples are the first ones containing {\it Gaussian geodesics}, and for the first time we prove embeddedness of the surfaces CSSCFF and CSSCCC from Callahan-Hoffman-Meeks-Wohlgemuth. 
\end{abstract}
\ \\
{\bf 1. Introduction}
\\

For a complete embedded minimal surface $S$ in $\R^3$ with finite total curvature, after the works from Rick Schoen \cite{RS} and Lopez-Ros \cite{LR} it has been known that examples with genus zero or number of ends $n\le 2$ are only possible for the plane and the catenoid. Therefore, new such surfaces should have at least genus one and three ends. Such a first example was found by Costa \cite{C}, followed by Hoffman-Meeks \cite{HM} still with three ends but then arbitrary genus. Moreover, in \cite{HM} the authors launched the conjecture that for all such $S$ it holds $n\le$ genus $+2$, which remains open until nowadays.
\\

In 1989, Karcher presented several examples in \cite{Ka1} and \cite{Ka2} that answered many questions in the theory of minimal surfaces. For instance, he presented the first such surfaces with positive genera and helicoidal ends, proved the existence of Alan Schoen's \cite{AS} triply periodic surfaces and gave doubly as well as singly periodic examples out of Scherk's families. By the way, after taking the quotient by the translation group, for saddle towers he obtained examples with $n=2k$ ends for genera zero and one, where $k\in\N$ and $k\ge$ genus $+2$. It is curious that no other {\it explicit} saddle towers were obtained since then, except for \cite{MRB}. This could be due to restrictions on such surfaces. For instance, Meeks and Wolf recently proved in \cite{MW} that a properly embedded saddle tower with four ends must belong to Scherk's family.
\\

In this paper, by an {\it almost explicit} example  we mean that one has the Weierstra\ss \  data. In addition, if the parameter's domain can always be refined for more and more precision, we call it an {\it explicit} example. In this sense, all Karcher's constructions are explicit. He structured them by a {\it reverse method}, which has a large literature of its application: \cite{BRB}, \cite{Ka1}, \cite{Ka2}, \cite{L}, \cite{LRB}, \cite{MRB}, \cite{V1}, \cite{V2}, \cite{V3}, \cite{V4}, \cite{V5} and \cite{V6}. Some non-explicit constructions are \cite{HMM}, \cite{Kp}, \cite{T3}, \cite{T2}, \cite{T1} and \cite{W2}. 
\\

Herewith we present the first saddle towers of genus two and $8=2\cdot 2+4$ ends. This might lead us to think about a sort of Hoffman-Meeks' conjecture for the slab, namely $n\ge 2($genus$+2)$. However, \cite{HMM} might throw some light upon Hoffman-Meeks' conjecture to answer it in the {\it negative}, because there the authors present saddle towers with arbitrary genus and three ends.  
\\
\input epsf
\begin{figure} [ht]
\centerline{
\hspace{-2cm}
\epsfxsize 8cm
\epsfbox{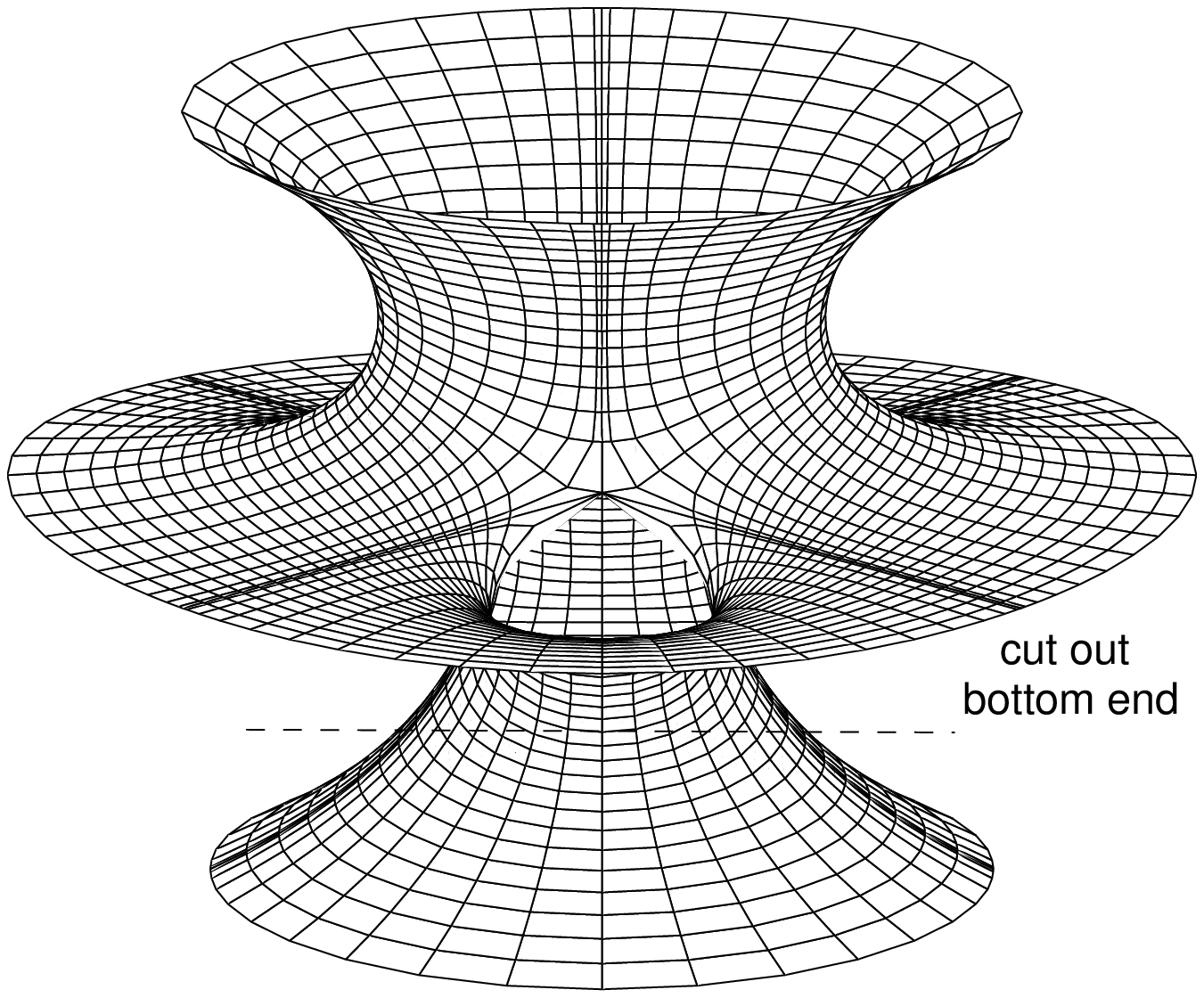}}
\hspace{8.5cm}(a)

\centerline{
\epsfxsize 8cm
\epsfbox{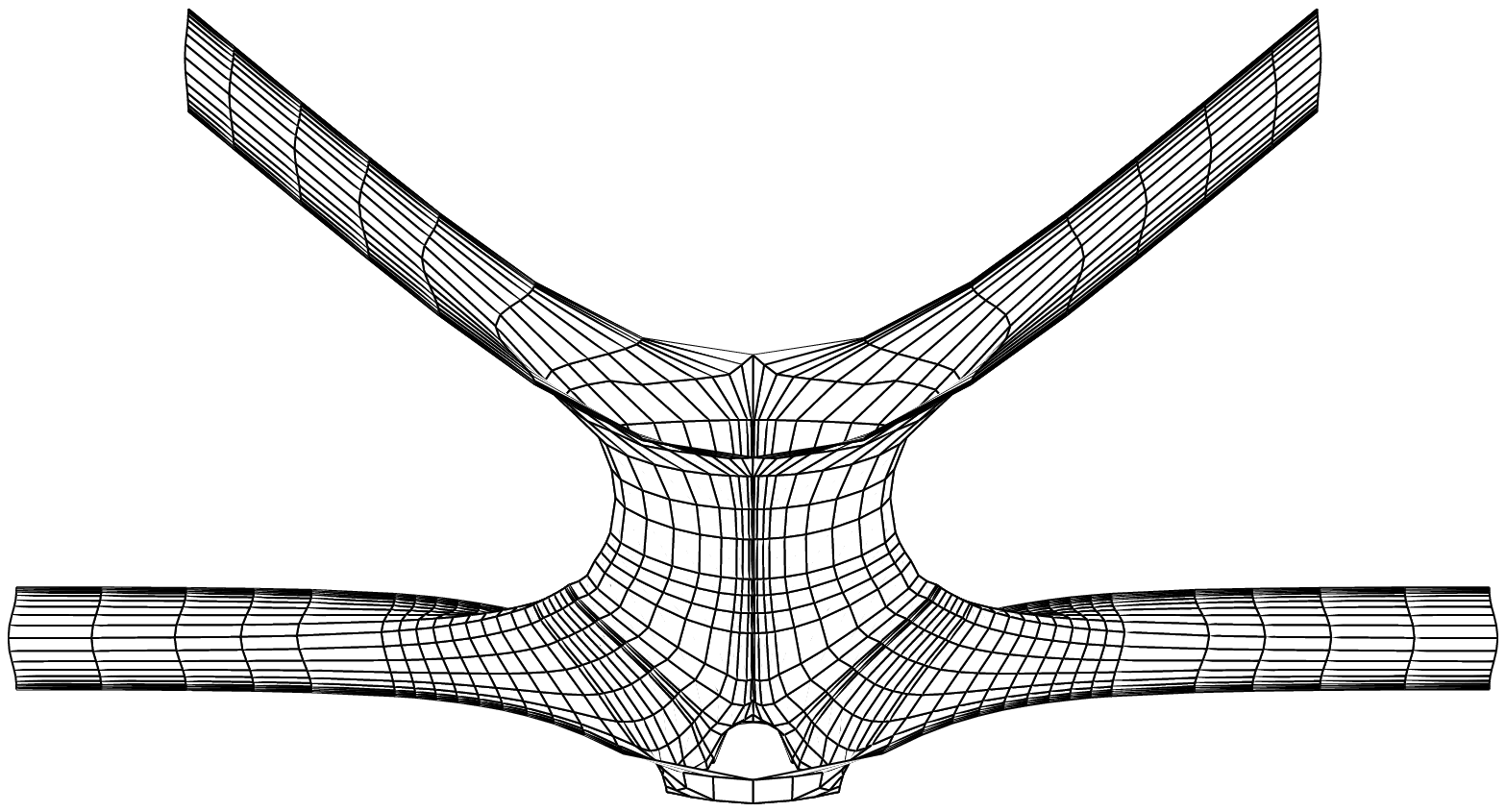}
\epsfxsize 9cm
\epsfbox{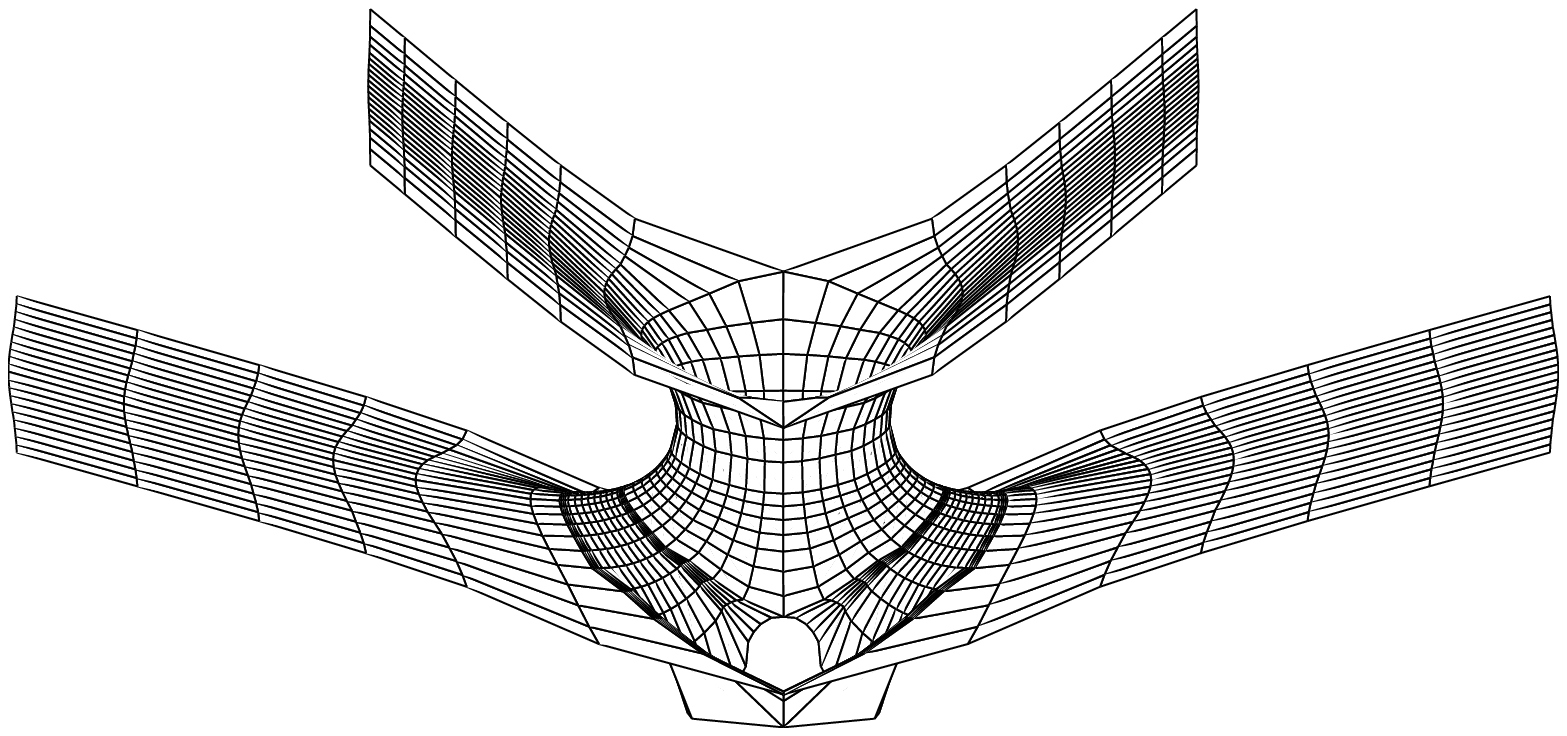}}
\hspace{4.5cm}(b)\hspace{8cm}(c)
\caption{(a) The Costa surface; (b) a special piece of $S$; (c) a general piece of $S$.}
\end{figure}

Our surfaces are easy to understand from Figures 1 and 2. Take the Costa surface and cut out its bottom catenoidal end, replacing it by a closed curve of reflectional symmetry. Afterwards, replace the remaining ends by Scherk-ends, as shown in Figure 1(b).
\\

Figure 2 represents the saddle towers we are going to construct. After we get the Weierstra\ss \  data by Karcher's method, there will be three period problems to solve, and this will follow practically {\it without computations}. This because we shall apply the {\it limit-method} described in either \cite{L} or \cite{LRB}. Other methods that ease the handling of period problems are found in \cite{BRB}, \cite{MRB} and \cite{W1}.     
\\
\input epsf
\begin{figure} [ht]
\epsfxsize 12cm
\epsfbox{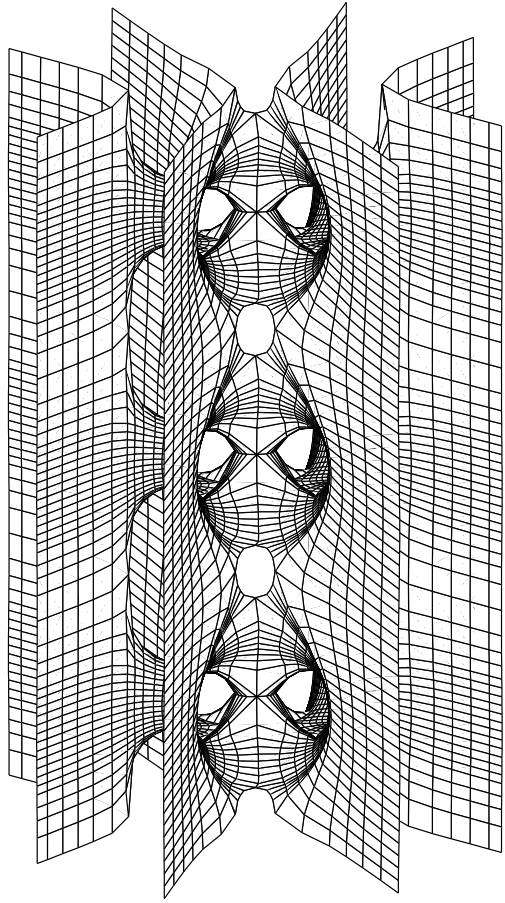}
\caption{A saddle tower of genus two.}
\end{figure}

In Figures 1(b) or 2, one notices the presence of a {\it Gaussian geodesic}. By this concept we mean a planar curve of reflectional symmetry, which is the graph of an even real-analytic function $f:\R\to[-1,0)$, where $f(0)=-1$, $f'\ne 0$ in $\R^*$ and $\Lim{t\to\infty}{f(t)}=0$. Since 1997, when the second author started his doctoral studies in Germany, he observed that they failed all construction attempts of minimal surfaces containing a Gaussian geodesic. In total one tried 15 different examples and periods never closed.
\\

This fact is important because, before our present work, for embedded minimal surfaces, it has been observed that the shape of an {\it unbounded} planar geodesic {\it always} matched one of the first nine standards in Figure 3, and {\it only} those. For each standard the picture cites an example, but the last one was missing. In fact, such geodesics seem to have a very restrictive geometry, and its study may reveal a lot about the general behaviour of the minimal surfaces. 
\\
\input epsf
\begin{figure} [ht]
\epsfxsize 12cm
\centerline{\epsfbox{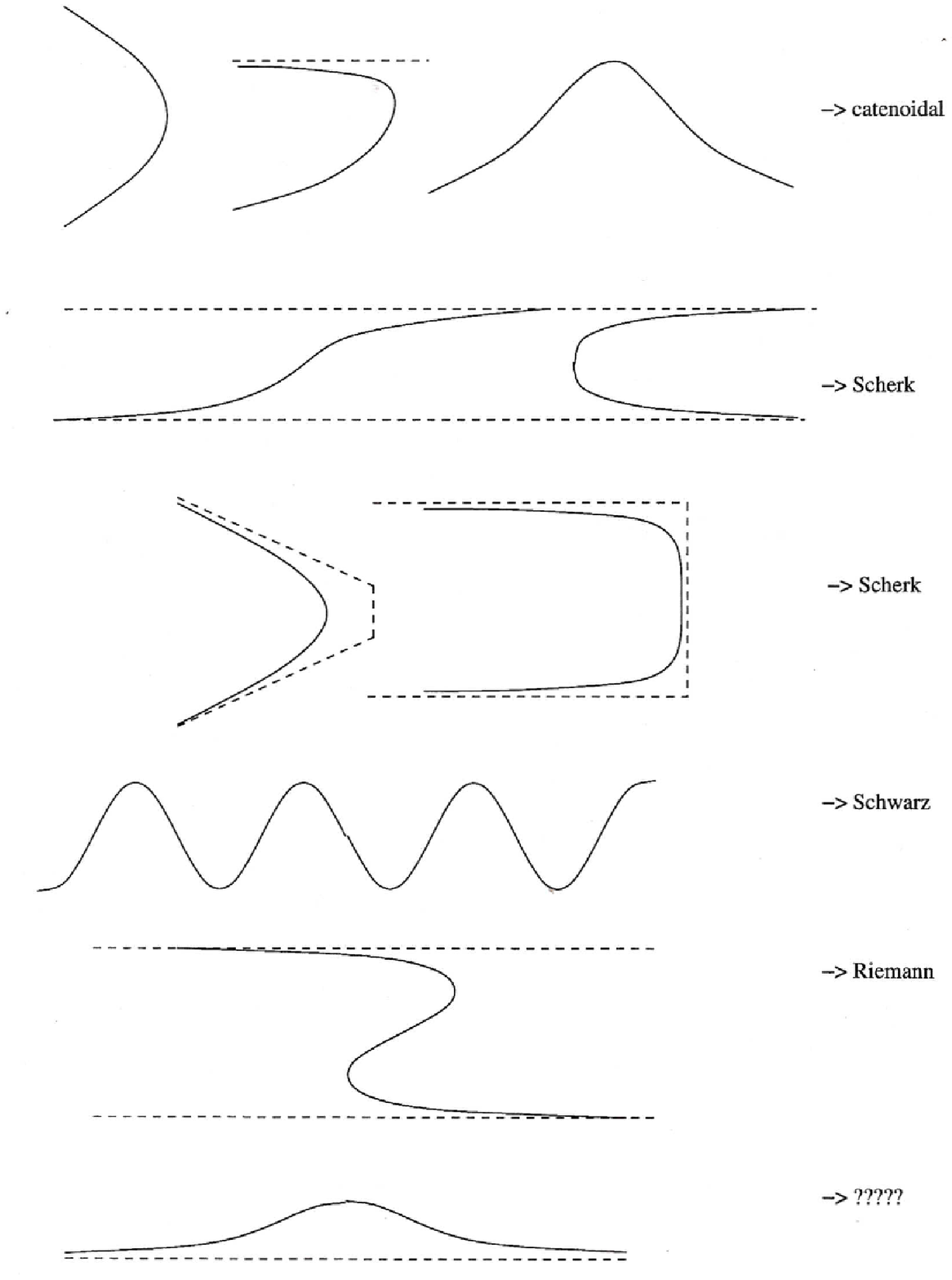}}
\caption{Reflectional symmetry curves on embedded minimal surfaces.}
\end{figure}

Our present paper then answers an open question cited in the introduction of \cite{V6}. Moreover, for the first time we also prove that the surfaces CSSCFF and CSSCCC, described in \cite{Wo}, are {\it embedded}. They will be used here as {\it limit-surfaces} for the method explained in either \cite{L} or \cite{LRB}. 
\\

Now we present the main theorem of this paper:
\\

{\bf Theorem 1.1.} \it There exists a continuous two-parameter family of saddle towers in $\R^3$, of which any member has the following properties:

i) The quotient by its translation group has genus two and eight Scherk-ends;

ii) It is invariant under reflections in $Ox_2x_3$, $Ox_3x_1$ and $mT/2+Ox_1x_2$, where $m\in\Z$ and $T$ is the single period of the surface.

iii) It is embedded in $\R^3$.
\\
Moreover, the family contains a continuous one-parameter sub-family of which any member has Gaussian geodesics.\rm
\\

In this present work, the second author was supported by the grants ``Bolsa de Produtividade Cient\ih fica'' from CNPq - Conselho Nacional de Desenvolvimento Cient\ih fico e Tecnol\'ogico, and FAPESP 04/02038-6.
\\
\\
{\bf 2. Preliminaries}
\\

In this section we state some basic definitions and theorems. Throughout this work, surfaces are considered connected and regular. Details can be found in \cite{JM}, \cite{Ka2}, \cite{LM}, \cite{N} and \cite{O}.
\\

{\bf Theorem 2.1.} \it Let $\X:R\to\E$ be a complete isometric immersion of a Riemannian surface $R$ into a three-dimensional complete flat space $\E$. If $\X$ is minimal and the total Gaussian curvature $\int_R K dA$ is finite, then there exists a compact Riemann surface $\ovl{R}$ and a finite number of points $p_1,\dots,p_r$ such that $R$ and $\ovl{R}\setminus\{p_1,\dots,p_r\}$ are conformally diffeomorphic.\rm 
\\

{\bf Theorem 2.2.} (Weierstra\ss \  representation). \it Let $R$ be a Riemann surface, $g$ and $dh$ meromorphic function and 1-differential form on $R$, such that the zeros of $dh$ coincide with the poles and zeros of $g$. Suppose that $\X:R\to\E$, given by
\BE
   \X(p):=Re\int^p(\phi_1,\phi_2,\phi_3),\eh\eh where\eh\eh
   (\phi_1,\phi_2,\phi_3):=\m\left(1/g-g,i/g+ig,2\right)dh,\label{W_rep}
\EE
is well-defined. Then $\X$ is a conformal minimal immersion. Conversely, every conformal minimal immersion $\X:R\to\E$ can be expressed as (1) for some meromorphic function $g$ and 1-form $dh$.\rm
\\

{\bf Definition 2.1.} The pair $(g,dh)$ is the \it Weierstra\ss \  data \rm and $\phi_1$, $\phi_2$, $\phi_3$ are the \it Weierstra\ss \  forms \rm on $R$ of the minimal immersion $\X:R\to\X(R)\subset\E$.
\\

{\bf Definition 2.2.} A complete, orientable minimal surface $S$ is \it algebraic \rm if it admits a Weierstra\ss \ representation such that $R=\ovl{R}\setminus\{p_1,\dots,p_r\}$, were $\ovl{R}$ is compact, and both $g$ and $dh$ extend meromorphically to $\ovl{R}$.
\\

{\bf Definition 2.3.} An {\it end} of $S$ is the image of a punctured neighbourhood $V_p$ of a point $p\in\{p_1,\dots,p_r\}$ such that $(\{p_1,\dots,p_r\}\setminus\{p\})\cap\ovl{V}_p=\emptyset$. The end is \it embedded \rm if this image is embedded for a sufficiently small neighbourhood of $p$.
\\

{\bf Theorem 2.3.} \it Under the hypotheses of Theorems 2.1 and 2.2, the Weierstra\ss \ data $(g,dh)$ extend meromorphically on $\ovl{R}$.\rm
\\

{\bf Theorem 2.4.} \it If in some holomorphic coordinates of a minimal immersion $F:\Omega\to\mathbb{R}^3$ there is a curve $\af$ such that the Gau\ss \ image $g\circ\af$ is contained either in a meridian or in the equator of $S^2$, and if also $dh(\dot{\af})$ is contained in a meridian of $S^2$, then $F\circ\af=\gamma$ is either in a plane or in a straight line (and is therefore a geodesic in both cases). The first case occurs exactly when $dh\cdot dg/g\in\R$ and the second when $dh\cdot dg/g\in i\R$. By the Schwarz Reflection Principle, if $\gamma$ is a straight line segment of $F(\Omega)=S$, then $S$ is invariant by $180^\circ$-rotation around $\gamma$. If $\gamma$ is a planar geodesic (not straight), then $S$ is invariant by a reflection in the plane of $\gamma$.\rm  
\\

{\bf Theorem 2.5.} (Jorge-Meeks Formula \cite{JM}). \it Under the hypotheses of Theorem 2.1, if the ends of $\X(R)=S$ are embedded, then
\[
deg(g)=k+r-1,
\]
where $k$ is the genus of $\ovl{R}$, and $r$ is the number of ends of the surface $S$.\rm  
\\
\\
REMARK: From the demonstration of Theorem 2.5, for the case of {\it Scherk-ends} the variable $r$ counts them in {\it pairs}. The function $g$ is the stereographic projection of the Gau\ss \ map $N:R\to S^2$ of the minimal immersion $\X$. It is a branched covering map of $\hat\C$ and $\int_SKdA=-4\pi$deg$(g)$. These facts will be largely used throughout this work.
\\
\\
{\bf 3. The compact Riemann surfaces $\ovl{M}$ and the Weierstra\ss \ data}
\\

Following Karcher's method, we are going to read off the necessary conditions for $S$ to exist. Afterwards we prove that these conditions are also sufficient. Figure 4 represents half of a {\it fundamental piece} $P$ of $S$, namely one that generates the whole surface by applying its translation group.
\\
\input epsf
\begin{figure} [ht]
\centerline{
\epsfxsize 16cm
\epsfbox{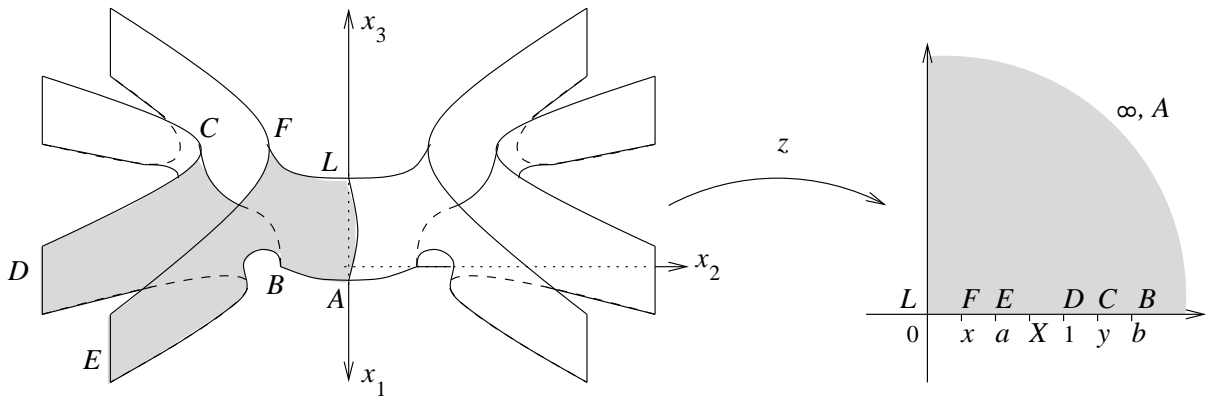}}
\caption{Half of a fundamental piece of $S$ and the $z$-map.}
\end{figure} 

A compactification of the Scherk-ends turns $P$ into a compact Riemann surface $\ovl{M}$ of genus two. Its {\it hyperelliptic involution} can be viewed as a 180$^\circ$-rotation about $Ox_2$. We call $\rho$ the quotient map induced by this rotation. Since $\ovl{M}/\rho$ is topologically $S^2$, Koebe's theorem together with a suitable M\"obius transformation will give a meromorphic function $z:\ovl{M}\to\hC$ such that $z(A)=\infty$, $z(D)=1$ and $z(L)=0$. The branch points of $z$ are $B$, $C$, $F$ and its images by reflection in $Ox_1x_3$. 
\\

Under the hyperelliptic involution, the reflection in either $Ox_1x_2$ or $Ox_2x_3$ leads to the same anti-holomorphic involution in $\hC$, where it must be the conjugate of a M\"obius transformation. Its fixed-point set is therefore a circumference passing through 0, 1 and $\infty$, namely $\R\cup\{\infty\}$. Thus, except for $AL$, $z$ takes all stretches in Figure 4 to real intervals. Hence $z(B)=b$, $z(C)=y$, $z(E)=a$ and $z(F)=x$, where these real free parameters satisfy the inequalities $0<x<a<1<y<b$. The following equation gives an algebraic description of $\ovl{M}$:
\BE
   w^2=\frac{b+z}{b-z}\cdot\frac{x-z}{x+z}\cdot\frac{y+z}{y-z}.
\EE

Under the hyperelliptic involution, the reflection in $Ox_1x_3$ is again the conjugate of a M\"obius transformation. Its fixed-point set is a circumference passing through 0 and $\infty$, but orthogonal to the real axis. Hence, $z(AL)=i\R_+$, and {\it not} \ $i\R_-$, because $z$ preserves orientation. Moreover, there is $X\in[a,1)$ such that the unitary normal $N$ is parallel to $Ox_2$ at $z^{-1}(X)$. The choice $X=a$ will give the Gaussian geodesics mentioned at the Introduction, whereas $X<a$ gives no embedded surfaces. 
\\

We choose the orientation of $S$ in order to have $g(A)=1$. Since $g$ is the stereographic projection of $N$, a careful analysis of the $z$- and $w$-divisors give $w(g+i)/(g-i)=(X+z)/(X-z)$, which implies 
\BE
   \biggl(\frac{g+i}{g-i}\biggl)^2=\frac{(b-z)(x+z)(y-z)(X+z)^2}{(b+z)(x-z)(y+z)(X-z)^2}.
\EE

We now list three possibilities that could theoretically occur:

(a) $N$ is vertical at some point of $BC$;

(b) $N$ is vertical at some point of $AL\setminus\{L\}$. 

(c) On $z(t)=t$, $b<t<\infty$, we have $g=e^{i\theta}$ with $\theta$ assuming non-negative values.
\\

Numerically, none of them happens. Anyway, the proof of Theorem 1.1 will follow independently. Regarding (a) and (b), numerical evidences make us expect that, for a certain complex $\af$ in the first open quadrant of $\C$, we have $g\in\{0,\infty\}$ only if $z=0$ or $z\in\{\pm\af,\pm\ovl{\af}\}$. In any case, 
\[
   0=(z-b)(z+x)(z-y)(z+X)^2+(z+b)(z-x)(z+y)(z-X)^2=2z(z^4+S_2z^2+S_4),
\]

where $S_2=by-bx-xy+2(x-b-y)X+X^2$ and $S_4=2bxyX+(by-bx-xy)X^2$. We also define $S_1=x+2X-b-y$, $S_3=bxy+2(by-bx-xy)X+(x-b-y)X^2$ and $S_5=bxyX^2$. In order to have $\pm\af$ and $\pm\ovl{\af}$ as roots, the following condition must hold: $S_2^2<4S_4$. However, such a restriction is not necessary to prove Theorem 1.1, as explained right above. Now one easily writes down the differential $dh$ as 
\BE
   dh=\frac{z^4+S_2z^2+S_4}{(z-b)(z-y)(z+x)}\cdot\frac{z\,dz/w}{(z^2-1)(z^2-a^2)}.
\EE

The following table summarises the behaviour of $g$ and $dh$ along important paths on $\ovl{M}$:
\[
\begin{tabular}{|c|c|c|c|c|}\hline 
Symmetry & Involution                   & $z$         & $g$     & $dh$  \\\hline
$AL$    & $(w,z)\to(1/\bar{w},-\bar{z})$& $i\R_+$     & $\R_+$  & $\R_+$ \\\hline
$LF$    & $(w,z)\to(\bar{w},\bar{z})$   & $0<\cdot<x$ & $i\R_+$ & $\R_+$ \\\hline
$FE$    & $(w,z)\to(-\bar{w},\bar{z})$  & $x<\cdot<a$ & $S^1$   & $i\R$ \\\hline
$ED$    & $(w,z)\to(-\bar{w},\bar{z})$  & $a<\cdot<1$ & $S^1$   & $i\R$ \\\hline
$DC$    & $(w,z)\to(-\bar{w},\bar{z})$  & $1<\cdot<y$ & $S^1$   & $i\R$ \\\hline
$CB$    & $(w,z)\to(\bar{w},\bar{z})$   & $y<\cdot<b$ & $i\R_-$ & $\R_-$ \\\hline
$BA$    & $(w,z)\to(-\bar{w},\bar{z})$  & $b<\cdot$   & $S^1$   & $i\R$ \\\hline
\end{tabular}
\]

From Theorem 2.4, we know that our surfaces {\it do} have the sought after symmetry curves to prove item (ii) of Theorem 1.1. However, we still need to solve the period problems to prove both items (i) and (ii). This is done in the next section. 
\\
\\
{\bf 4. The period problems}
\\

From Figure 4 and Section 3 we may write down the residue and period problems. An easy computation gives:
\BE
   r:=2\pi i\,Res(dh,D)=\frac{\pi(1+S_2+S_4)}{(1-a^2)\sqrt{(1-b^2)(1-x^2)(1-y^2)}},
\EE
\BE
   R:=2\pi i\,Res(dh,E)=\frac{\pi(a^4+S_2a^2+S_4)}{(1-a^2)\sqrt{(a^2-b^2)(a^2-x^2)(a^2-y^2)}}.
\EE

The residue problem will be solved if both (5) and (6) match. We recall (c) and the fact that, at least numerically, it does not happen to the limit-surfaces CSSCFF and CSSCCC. Therefore, we expect that $-\pi/4<\theta<0$. Together with (3) and (4), we compute 
\BE
   I:=\m\int_{BA}(1/g-g)dh=\int_b^\infty\frac{(S_1t^4+S_3t^2+S_5)dt}{(t^2-1)(t^2-a^2)\sqrt{(t^2-b^2)(t^2-x^2)(t^2-y^2)}}.
\EE

For $z(t)=it$, $0<t<\infty$, we always have $g>0$, but {\it not} always $g>1$. This is what happens to the limit-surfaces, at least numerically. Define
\BE 
   J:=\m\int_{LA}(1/g-g)dh=\m\int_0^\infty\frac{t(t^4-S_2t^2+S_4)(g-1/g)dt}{(t^2+1)(t^2+a^2)\sqrt{(t^2+b^2)(t^2+x^2)(t^2+y^2)}}.
\EE

The first period problem will be solved providing $I=J$. At last, define
\BE
   K:=\int_{BC}dh=\int_y^b\frac{t^4+S_2t^2+S_4}{(1-t^2)(a^2-t^2)}\cdot\frac{t\,dt}{\sqrt{(b^2-t^2)(x^2-t^2)(y^2-t^2)}}.
\EE 
The second period problem will be solved if $K=R/2$.    
\\
\\
{\bf 5. Solution of the period problems}
\\

For $(a,b,x,X,y)\in\R^5$, we consider the function
\BE
   F:=(1+S_2+S_4)^2(a^2-x^2)(a^2-b^2)(a^2-y^2)-(a^4+S_2a^2+S_4)^2(1-y^2)(1-b^2)(1-x^2).
\EE

An easy computation gives
\BE
   \frac{\deh F}{\deh y}\biggl|_{(x,X,y)=(a,a,1)}=2(a^4+S_2a^2+S_4)^2(1-b^2)(1-a^2).
\EE

It is immediate to see that $a^4+S_2a^2+S_4=4a^2(1-a)(b-a)\ne 0$. Therefore, by the {\it implicit function theorem}, there is a unique function $y=y(a,b,x,X)$ that makes $F\equiv 0$ for $(a,b)\in(0,1)\times(1,+\infty)$ and $a-\eps<x,X<a+\eps$, for a certain $\eps=\eps(a,b)>0$. Now restrict the variables $(a,b)$ to $(0,1)\times(1,+\infty)$, and the variables $(x,X,y)$ to $(0,a)\times(a,1)\times(1,b)$. In view of (5) and (6), we shall have $r=R$ providing $F\equiv 0$ for the choice of $y$ as the {\it implicit function} of $(a,b,x,X)$. Of course, we are interested in the case $x<a\le X$, for which we should guarantee that $y>1$. But this comes directly from (10), since $1+S_2+S_4=(b+1)(1-a)^3>0$.  
\\

From this point on, we take $y$ as the function obtained above. Now define 
\BE
   G:=-i\cdot\frac{g+i}{g-i}\eh\eh{\rm and}\eh\eh dH:=\frac{i}{2}(1/g+g)dh.
\EE

In terms of a rigid motion in $\R^3$, the Weierstra\ss \ data $(G,dH)$ provide the {\it same} minimal surfaces from $(g,dh)$, but rotated counterclockwise by 90$^\circ$ around $Ox_1$. Indeed, one easily checks that
\[
   (1/G-G)dH=(1/g-g)dh\eh\eh{\rm and}\eh\eh i(1/G+G)dH=-2\,dh.
\]

Moreover, a simple reckoning gives
\BE
   dH=\frac{(X^2-z^2)dz}{(a^2-z^2)(1-z^2)}.
\EE

From this point on we shall strongly use the reference \cite{Wo}. Take small and disjoint open neighbourhoods $U\ni a$ and $V\ni 1$. The set $\K=\ovl{M}\setminus z^{-1}((\pm U)\cup(\pm V))$ is then compact. In the case $X=a$, from (3), (12) and (13) one sees that $(G,dH)$ converges uniformly on $\K$ to the Weierstra\ss \ data of the surfaces CSSCFF, as described in \cite{Wo}, page 16. Figure 5 is a reproduction of the same picture in \cite{Wo}, page 16, but conveniently placed and marked with points in order to visualise what happens to our surfaces for the extreme values $x=a$ and $y=1$. Notice that we {\it have} fixed $X=a$.
\\
\input epsf
\begin{figure} [ht]
\centerline{
\epsfxsize 8cm
\epsfbox{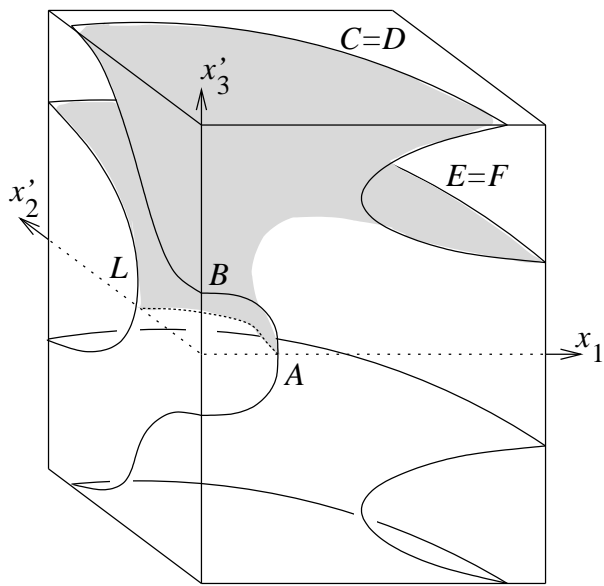}}
\caption{One quarter of the CSSCFF surface with $x_2^\prime=x_3$ and $x_3^\prime=-x_2$.}
\end{figure}

In \cite{Wo}, page 17, the author defines the periods $\pi_1(a,b)$ and $\pi_2(a,b)$. The first one is the integral of $\phi_1$ along an upper arc connecting some point in $(1,b)$ to some point in $(0,a)$. The second is the integral of $\phi_2$ along an upper arc connecting some point in $(b,\infty)$ to some point in $(a,1)$. The first arc is homotopically $BAL$, and the second is homotopically the oriented real segment from $b$ to $X$, which gives a Cauchy-principal value for the integral of $\phi_2$. However, since the integrals are invariant by free homotopy, after taking the extreme values $x=X=a$ and $y=1$, the integrals $I-J$ and $K-R/2$ will match $\pi_1$ and $\pi_2$, respectively. The following table summarises his conclusions about $\pi_{1,2}$ for $(a,b)\in[0,1]\times[1,\infty]$:
\[
\begin{tabular}{|c|c|c|c|c|}\hline 
$a$         & $b$       & $\pi_1$   & $\pi_2$  \\\hline
$0<\cdot<1$ & $1$       & chs.sign  & $<0$     \\\hline
$1$         & $1<\cdot$ & $<0$      & chs.sign \\\hline
$1>\cdot>0$ & $\infty$  & $>0$      & $>0$     \\\hline
$0$         & $\cdot>1$ & $>0$      & $>0$     \\\hline
$0$         & $1$       & $+\infty$ & $0$      \\\hline
$1$         & $1$       & $-\infty$ & $<0$      \\\hline
$1$         & $\infty$  & $0$       & $+\infty$ \\\hline
$0$         & $\infty$  & $1$       & $1$       \\\hline
\end{tabular}
\]

Moreover, he proves that the graphs of $\pi_{1,2}$ intersect along a space curve $\cC_0$, which in its turn has a crossing with the level-zero horizontal plane. This crossing occurs at a point $(a_0,b_0)$ and solves the period problems for CSSCFF. Such a fact is numerically represented by Figure 6, where $0.82\le a\le 0.98$ and $1.01\le b\le 1.51$. We remark that the numeric solution occurs for $(a,b)$ quite close to $(1,1)$.
\\
\input epsf
\begin{figure} [ht]
\centerline{
\epsfxsize 7cm
\epsfbox{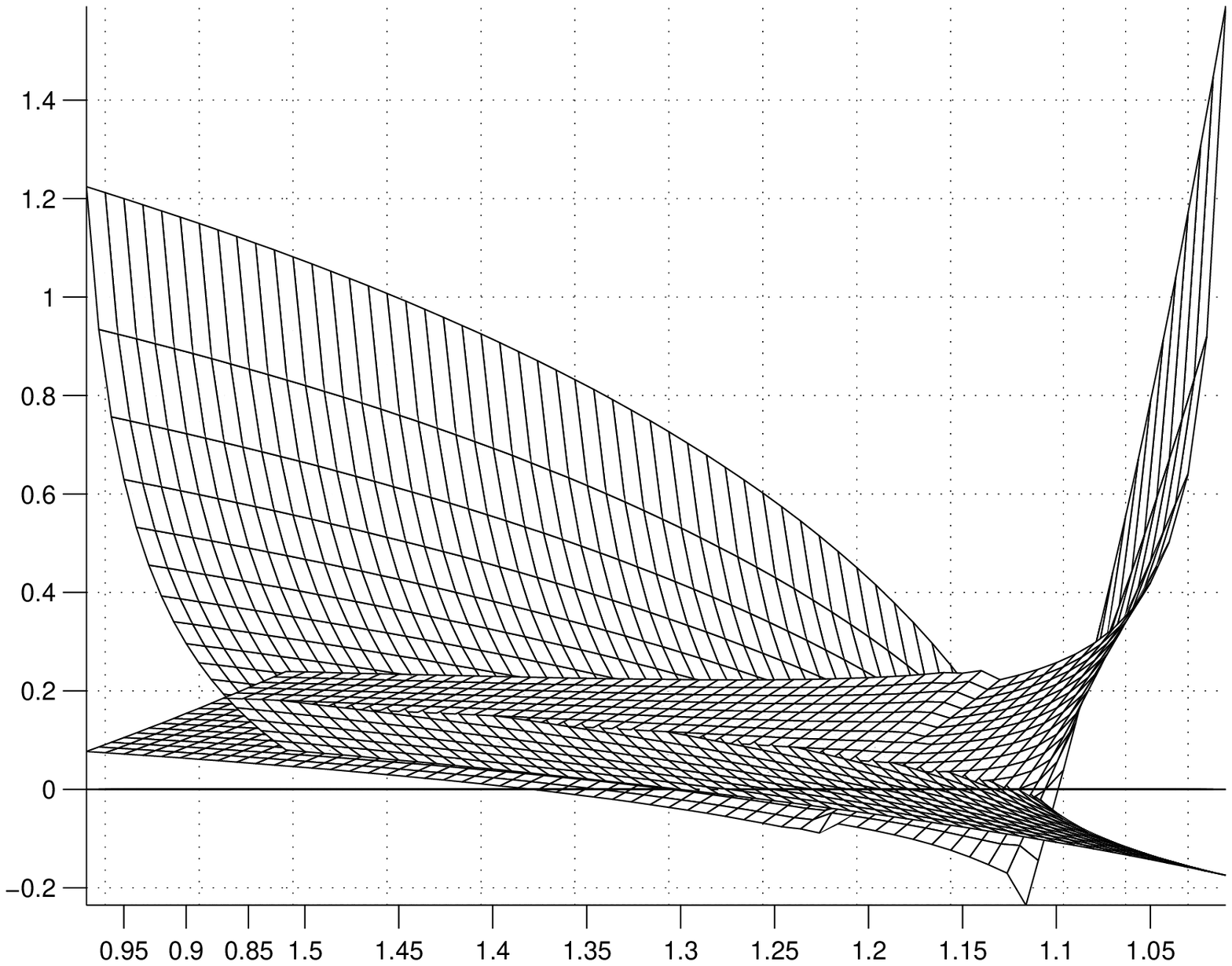}
\epsfxsize 7.6cm
\epsfbox{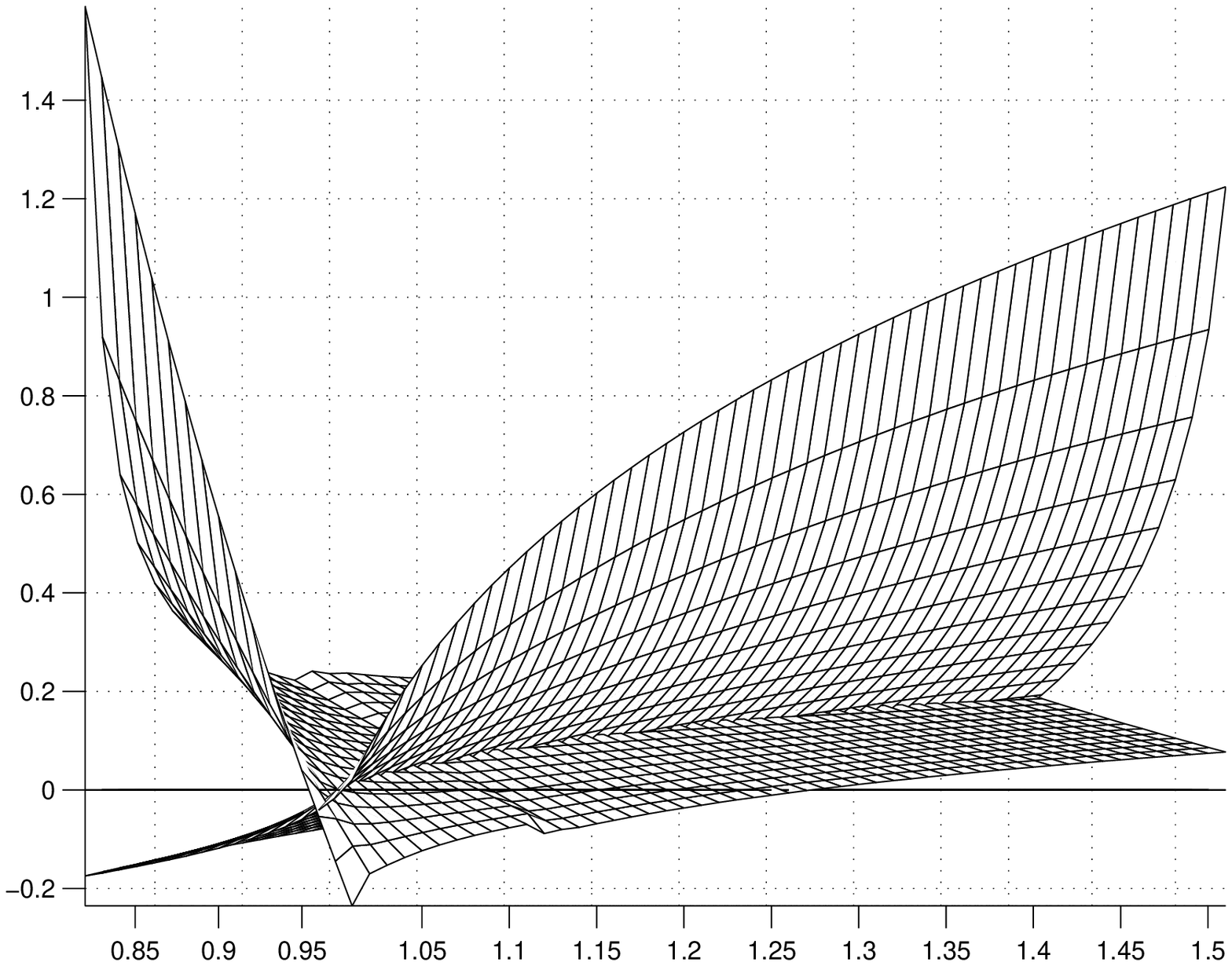}}
\caption{Two views of the $\pi_{1,2}$-graphs.}
\end{figure}

In Section 5, for each pair $(a,b)\in(0,1)\times(1,+\infty)$ we obtained a function $y(a,b,x,X)$ that makes $F\equiv 0$ for $a-\eps<x,X<a+\eps$ and a certain $\eps=\eps(a,b)>0$. By taking a small relatively compact neighbourhood of $(a_0,b_0)$, we may assume that $\eps$ does not depend on $(a,b)$ in this neighbourhood. We fix $X=a$ and extend the definition of $\pi_{1,2}$ for $a-\eps<x<a$ as $I-J$ and $K-R/2$, respectively. Since this extension is smooth, we shall get a continuous family of space curves $\cC_t$, each of them still crossing the horizontal plane at a point $(a_t,b_t)$, where $t\in[0,\eps)$. Each crossing happens for $x$ taken as the function $x(a)=a-t$, $a\in(0,1)$.       
\\

From \cite{Wo}, page 21, by fixing $X=\tilde{a}$ the same reasoning applies now with the surfaces CSSCCC. For $\eps$ small enough we get a continuous {\it two}-parameter family of saddle towers in $\R^3$, tracked by $(t,X)\in[0,\eps)\times[a,a+\eps)$, with the properties (i) and (ii) described in Theorem 1.1. The Gaussian geodesics occur for $X=a$, and this fact is easily checked by (3) along $z\in[x,a)$, since the curve is planar and its unitary tangent vector is just a clockwise rotation of $g$ by 90$^\circ$.
\\
\\
{\bf 6. Embeddedness}
\\

In the previous section, we proved that our surfaces are period-free in the slab. From Section 3, all the sought after symmetry curves were verified to exist. In particular, the behaviour of the Gau\ss \ map is summarised in Figure 7(a). In this picture, we have stressed the {\it inner} points where $|g|=1$. Moreover, notice the branches of $g$ when the normal passes along $LA$, $BC$ and from $z=X$ to $F$, as sub-stretch of $DEF$. The shaded region, which bounds $g(i\R_+\cup\R_+)$, is determined by the fact that deg$(g)=5$, according to Theorem 2.5. Figure 7(a) depicts this image under $g$. Moreover, with exactly 8 copies of Figure 7(a) we cover $\hC$ five times, of course by taking plus-minus conjugates and inversion with respect to $S^1$. This confirms the correct choice of the shaded region. Indeed, for if it were in the complement, we would get a contradiction with deg$(g)=5$. 
\\
\input epsf
\begin{figure} [ht]
\centerline{
\epsfxsize 16cm
\epsfbox{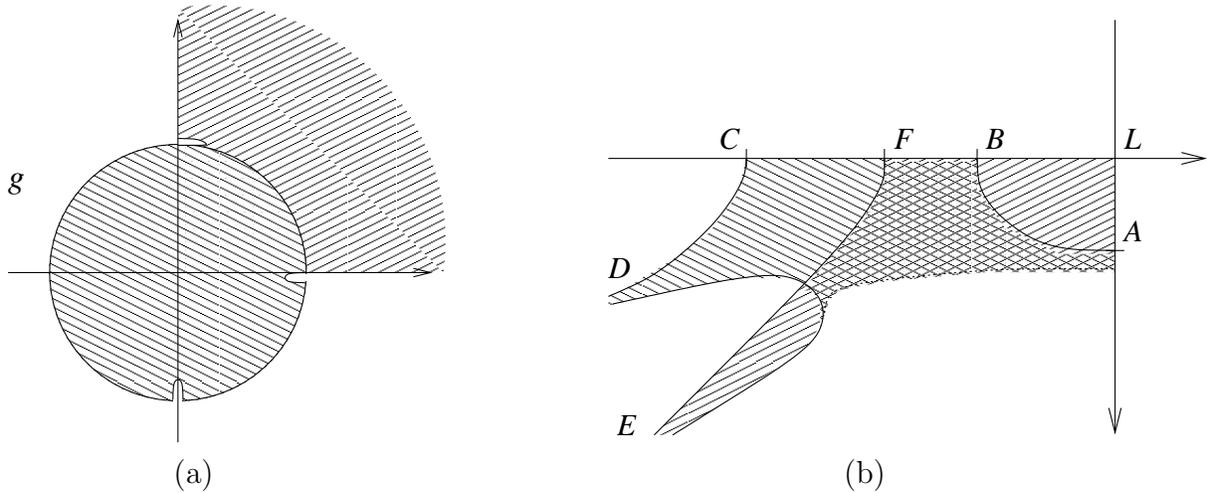}}
\hspace{1.2in}(a)\hspace{3.3in}(b)
\caption{a) $g$-image of $z$ in the first quadrant; b) corresponding projection onto $x_3=0$.}
\end{figure}

Of course, the branch on $BC$ could extend till pass by $g=0$, and the one from $LA$ could come out on $AB$ instead. In this case, it would cling to $S^1$, and not to $\R$ as shown in Figure 7(a). These possibilities correspond to (a) and (c), listed in Section 3. Regarding (b), it would extend the branch on $LA$ till it pass by $g=0$. Independently of these particular possibilities, embeddedness will follow anyway.
\\

Regarding Figure 7(b), it is the {\it expected} projection of regions $|g|\ge 1$ and $|g|\le 1$ onto $x_3=0$, and they overlap in the darkest shaded sub-region. Some other possibilities are shown in Figure 8, and we could even add cases in which $B$ comes out between $C$ and $F$, or even $B=F$. Its position, however, will not be crucial to our demonstration, although $CFB$ is numerically correct. Except for this particular detail, we shall prove that Figure 7(b) is the only possible one. 
\\
\input epsf
\begin{figure} [ht]
\centerline{
\epsfxsize 16cm
\epsfbox{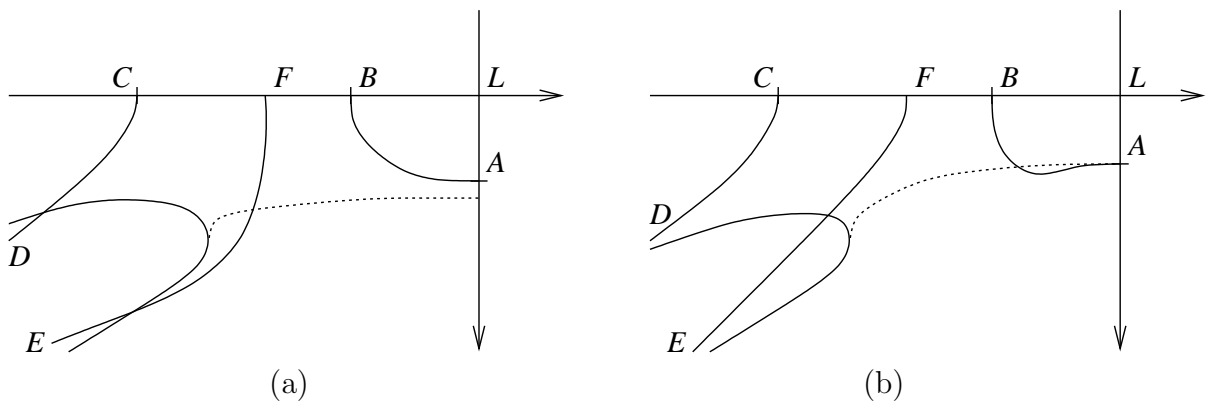}}
\hspace{1.7in}(a)\hspace{2.9in}(b)
\caption{Variations of Figure 7(b).}
\end{figure}

Among the situations (a)-(c), listed in Section 3, whether any of them occurs, the fact is that $g$ has a branch somewhere along $LAB$, another in $BC$, and finally a third one in $(x,X)\ni z$. On $\ovl{M}$, this gives a {\it total} of 12 zeros for $dg$, because deg$(dg)=2$, and deg$(g)=5$ implies 10 poles for $dg$. This means, $g$ is unbranched inside the shaded region in Figure 7(a), which does not count the contour $g(i\R_+\cup\R_+)$.
\\

We recall that exactly 8 copies of Figure 7(a) cover $\hC$ five times. This means, $g$ is {\it injective} for $z$ in the first {\it open} quadrant of $\hC$. Therefore, there is a simple curve $\Gamma$ in the first open quadrant of which the extremes are $z=X$ and a certain $Y\in(b,\infty]\cup i\R_+^*$, such that $g(\Gamma)$ is the inner unitary arc in Figure 7(a). The curve $\Gamma$ divides that quadrant in two disjoint components $\A$ and $\B$, corresponding to $|g|<1$ and $|g|>1$ in Figure 7(a), respectively. 
\\

Since $g$ is injective in the first open quadrant, from our minimal immersion $\X=(x_1,x_2,x_3)$ in the slab, given by Theorem 2.2, we have that $(x_1,x_2):\A\to\R^2$ and $(x_1,x_2):\B\to\R^2$ are {\it both immersions}. Moreover, their images are {\it connected} open subsets of $\R^2$. Therefore, any two paths of the $\deh\A$-image under $(x_1,x_2)$  are {\it disjoint}, and the same holds for $\deh\B$. Because of that, among the three intersections depicted in Figure 8(a), none of them occurs, not even as tangent points. 
\\

For the same reason, the dotted curve in Figure 8(b), which represents the image of $\Gamma$ under $(x_1,x_2)$, could not cross any of the continuous curves represented there. Instead, its right-hand side extreme should be tangent to the {\it uppermost} point of $AB$ (notice that $Ox_1$ is vertical {\it downwards} in the picture). Besides, Scherk-ends are asymptotic to half-planes, which are vertical in our case, so that the filled up regions {\it must} be the insides indicated by Figure 7(b).    
\\

Let us now analyse the convexity and monotonicity of the patches that build $g(i\R_+\cup\R_+)$. Henceforth, any patch will be viewed as its {\it projection} onto $x_3=0$, and we shall strongly use the fact that $g$ {\it is the stereographic projection of the Gau\ss \ map}. Since the limit-surfaces are CSSCFF and CSSCCC, take ours so that ends $E$ and $D$ do not intersect, and $F$ is at the right of $C$. Moreover, we {\it know} the branches of $g$. In particular, stretch $DE$ is a monotone convex curve. It is {\it not} always true that $FE$ is convex, but for sure monotone, because there we have $g=e^{i\theta}$, $-\pi/2<t\le 0$. By the way, we get the {\it Gaussian geodesics} exactly when $X=a$.
\\

The curve $g(\Gamma)$ is convex and monotone, as well as $CD$ and $FL$. From the cases listed in Section 3, it would fail monotonicity for $BC$ in case (a), {\it unless} we have a monkey-saddle where $N$ is vertical, and convexity for $AB$ in case (c) {\it if} $\theta$ assumes positive values. Anyway, they are {\it always} convex and monotone, respectively, and one may use the Weierstra\ss \ data, specially (13), to confirm these facts. If (c) does not occur, then $AL$ is not monotone {\it unless} we have a four-fold symmetry saddle at $A$. Anyway, it is {\it always} convex. 
\\

Let us now analyse $(x_1,x_2):\A\to\R^2$. By considering $\A\subset\hC$, we may continuously extend it to $(x_1,x_2):\ovl{\A}\to\hC$. The pre-image of any point in $(x_1,x_2)(\A)$ is a finite set of points, otherwise they would accumulate in some $p\in\deh\A$, a contradiction. Therefore, $(x_1,x_2)$ is a covering map from $\A$ to the simply connected region $(x_1,x_2)(\A)$. Namely, it is injective. By the same arguments, $(x_1,x_2):\B\to\R^2$ is also injective.
\\

From this point on, the patches will be viewed as space curves again. Since deg$(g)=5$, $LA$, $AB$ and $BC$ are each free of self-intersections. For the latters, we re-confirm this fact with (13). It might happen, however, that $BC\cap FL\ne\emptyset$. Now we recall that $(x_1,x_2)$ is an immersion for either $\A$ or $\B$, and therefore none of those sets could have an image point on $g(\Gamma)$. Consequently, except for the common stretch $\X(\Gamma)$, $BC\cap FL$ is the only possible intersection between the boundaries of the graphs $(\X|_\A)$ and $(\X|_\B)$.
\\

Therefore, $\X(\ovl{\A})\cap\X(\ovl{\B})\setminus\X(\Gamma)=(\X(\A)\cap\X(\B))\cup(BC\cap FL)$. Now we are going to apply the {\it maximum principle} (see \cite{RS}, for instance). If that set were not empty, then we could lift the graph $(\X|_\B)$ till getting a first contact point. The maximum principle would then imply that both pieces coincide, which is absurd. Therefore, our surface has an embedded {\it fundamental domain}, confined to a slab in the 4th octant of $\R^3$. By successive reflections in its boundary, we get an embedded singly periodic surface in $\R^3$. This last argument finally demonstrates Theorem 1.1.
\\

We conclude this section with the embeddedness of the surfaces CSSCFF and CSSCCC. In Section 5 we proved that our surfaces are parametrised by $t\in[0,\eps)$, where $t=0$ gives these limit-cases. In fact, any of them has disjoint embedded ends. So, if there were a self-intersection, this would happen inside a compact like $\K$ described in Section 5. Again by the {\it maximum principle} for minimal surfaces, the same would happen to our surfaces for some positive $t$ close enough to 0. But this is impossible because of Theorem 1.1. Therefore, the surfaces CSSCFF and CSSCCC are also embedded in $\R^3$.

\ \\
da Silva, M\'arcio Fabiano\\
Universidade Federal do ABC\\
r. Catequese 242, 3rd floor\\
09090-400 Santo Andr\'e - SP, Brazil\\
E-mail: {\tt marcio.silva@ufabc.edu.br}
\\
\\
Ramos Batista, Val\'erio\\
Universidade Federal do ABC\\
r. Catequese 242, 3rd floor\\
09090-400 Santo Andr\'e - SP, Brazil\\
E-mail: {\tt valerio.batista@ufabc.edu.br}
\end{document}